\documentclass[11pt]{amsart}
\usepackage{amscd}      % to be able to use "CD" environment
\usepackage{amsbsy, amsthm, eucal, graphics,epic,amsbsy, stmaryrd, epic, eepic}
\usepackage{amssymb}
\usepackage{amsmath}
\usepackage[all]{xy}
\usepackage[dvips]{color}

\usepackage{latexsym}
\usepackage{epsfig}
\usepackage{amsfonts}
\usepackage{enumerate}

\newtheorem{prop}{Proposition}[section]
\newtheorem{lem}[prop]{Lemma}

\newtheorem{cor}[prop]{Corollary}
\newtheorem{thm}[prop]{Theorem}
\newtheorem{rmk}[prop]{Remark}

\newtheorem{notation}[prop]{Notation}

\newtheorem{defn}[prop]{Definition}

\newenvironment{pf}{\begin{trivlist}\item[]{\sc Proof.}}%
            {\nolinebreak $\Box$ \end{trivlist}}

\newcommand{\noprint}[1]{}

\renewcommand{\tilde}{\widetilde}

\newcommand{\cl}{\mathop{\rm cl}\nolimits}

\newcommand{\id}{\mathop{\rm id}\nolimits}

\newcommand{\ldiag}[1]%
       {\makebox[0cm]{${\scriptstyle#1}\downarrow\phantom{\scriptstyle#1}$}}
\newcommand{\ldiagup}[1]%
       {\makebox[0cm]{${\scriptstyle#1}\uparrow\phantom{\scriptstyle#1}$}}
\newcommand{\rdiag}[1]%
       {\makebox[0cm]{$\phantom{\scriptstyle#1}\downarrow{\scriptstyle#1}$}}
\newcommand{\sediagr}[1]%
       {\makebox[0cm]{$\phantom{\scriptstyle#1}\searrow{\scriptstyle#1}$}}
\newcommand{\nediagr}[1]%
       {\makebox[0cm]{$\phantom{\scriptstyle#1}\nearrow{\scriptstyle#1}$}}
\newcommand{\rdiagup}[1]%
       {\makebox[0cm]{$\phantom{\scriptstyle#1}\uparrow{\scriptstyle#1}$}}
\newcommand{\swdiag}[1]%
       {\makebox[0cm]{$\phantom{\scriptstyle#1}\swarrow{\scriptstyle#1}$}}
\newcommand{\sediag}[1]%
       {\makebox[0cm]{${\scriptstyle#1}\searrow\phantom{\scriptstyle#1}$}}
\newcommand{\nediag}[1]%
       {\makebox[0cm]{${\scriptstyle#1}\nearrow\phantom{\scriptstyle#1}$}}

\newcommand{\doublearrowstack}[2]%
                      {{{{\scriptstyle#1}\atop{\textstyle\longrightarrow}}\atop{{\textstyle\longrightarrow}\atop{\scriptstyle#2}}}}
\newcommand{\rightleftarrowstack}[2]%
                      {{{{\scriptstyle#1}\atop{\textstyle\longrightarrow}}\atop{{\textstyle\longleftarrow}\atop{\scriptstyle#2}}}}
\newcommand{\leftrightarrowstack}[2]%
                      {{{{\scriptstyle#1}\atop{\textstyle\longleftarrow}}\atop{{\textstyle\longrightarrow}\atop{\scriptstyle#2}}}}

\newcommand{\overtoparrow}%
{\makebox[0cm]{\beginpicture \setcoordinatesystem units
<.8cm,.4cm> point at 0 0 \setplotarea x from -3 to 3, y from 0 to
1 \setquadratic \plot -3 0 0 1 3 0 / \put{\vector(3,-1){0}}[Bl] at
3 0
\endpicture}}

\newcommand{\underbottomarrow}%
{\makebox[0cm]{\beginpicture \setcoordinatesystem units
<.8cm,.4cm> point at 0 0 \setplotarea x from -3 to 3, y from 0 to
1 \setquadratic \plot -3 1 0 0 3 1 / \put{\vector(3,1){0}}[Bl] at
3 1
\endpicture}}

\newcommand{\ses}[5]%
{0\longrightarrow#1\stackrel{#2}{ \longrightarrow}#3\stackrel{#4}{
\longrightarrow}#5\longrightarrow0}

\newcommand{\dt}[6]%
{#1\stackrel{#2}{longrightarrow}#3
\stackrel{#4}{\longrightarrow}#5 \stackrel{#6}{\longrightarrow}
#1[1]}

\newcommand{\cat}[1]%
{(\mbox{\rm #1})}

  \def\bbC{{\mathbb C}}  

\def\bbG{{\mathbb G}} \def\bbP{{\mathbb P}}   
 \def\bbZ{{\mathbb Z}}   \def\bbN{{\mathbb N}}

\newcommand{\clA}{{\mathcal{A}}}
\newcommand{\clB}{{\mathcal{B}}}

\newcommand{\clE}{{\mathcal {E}}}

\newcommand{\clG}{{\mathcal{G}}}
\newcommand{\clF}{{\mathcal F}}

\newcommand{\clM}{\mathcal{M}}

\newcommand{\clH}{{\mathcal H}}

\newcommand{\clR}{{\mathcal R}}
\newcommand{\clT}{{\mathcal{T}}}
\newcommand{\clU}{{\mathcal{U}}}
\newcommand{\clO}{{\mathcal{O}}}
\newcommand{\clQ}{{\mathcal{Q}}}
\newcommand{\clX}{{\mathcal{X}}}
\newcommand{\clY}{{\mathcal{Y}}}
\newcommand{\clW}{{\mathcal{W}}}
\newcommand{\clZ}{{\mathcal{Z}}}
\def\longto{\longrightarrow}
\def\isomto{\stackrel{\sim}{\longto}}

\title{Moduli spaces of pairs over projective stacks}
\author{Elena Andreini}
\address{International School for Advanced Studies\\Via Beirut 2-4\\
34100 Trieste\\ Italy}
\email{andreini.elena@gmail.com}

\begin{document}
 %%% -----------------------------------------------------------------------
\bibliographystyle{alpha}

\sloppy \maketitle
%%% ----------------------------------------------------------------------
\tableofcontents

\begin{abstract}
Let $\clX$  a projective stack over an algebraically closed field $k$ of characteristic 0. Let $\clE$ be a generating sheaf over $\clX$ and $\clO_X(1)$ a polarization of its coarse moduli space $X$. We define a notion of pair 
which is the datum of a non vanishing morphism $\Gamma\otimes\clE\to \clF$ where $\Gamma$ is a finite dimensional $k$ vector space
and $\clF$ is a coherent sheaf over $\clX$. We construct the stack and the  moduli space of semistable pairs.
 The notion of semistability depends on a polynomial parameter and it is dictated by the GIT
construction of the moduli space. 
\end{abstract}
\section{Introduction}
Recently  a lot of attention has been drawn by sheaf theoretic curve counting theories of projective threefolds.
Among them Pandharipande-Thomas invariants \cite{PandThom1-09} are computed via integration over
the virtual fundamental class of the moduli space of the so called {\em stable pairs}.
The moduli spaces used to compute PT invariants are a special case of moduli spaces of coherent systems  introduced by Le Potier in \cite{LePot93}.  A coherent system is the datum  of a pair $(F, \Gamma)$, where $F$ is   a pure $d$-dimensional 
coherent sheaf  and $\Gamma\subset H^0(X,F)$  is  a subspace of its global sections. 
The moduli spaces are constructed as projective varieties via GIT techniques. The GIT stability
condition is equivalent to a modified   Gieseker stability, where the Hilbert polynomial is corrected 
by a contribution proportional to $\mbox{dim}\ \Gamma$ and to a polynomial stability parameter.  A coherent system  can be reconstructed from the associated evaluation morphism
$ev:\Gamma\otimes\clO_X\to F$. 
%By using this reformulation, semi-stable coherent systems are embedded in an abelian category with enough injetives. This approach  which is needed to study deformation theory. 
In this note we study a similar moduli problem over projective stacks. If we work over an algebraically closed field $k$ of 
characteristic zero, a projective stack is a stack with projective coarse moduli space that can be embedded into a smooth proper Deligne-Mumford stack. Any projective stack admits a {\em generating} sheaf $\clE$, namely
a locally free sheaf whose fibers carry every representation of the automorphism group of the underlying point.
We propose a notion of {\em pair} on projective stacks which is a natural generalization of the evaluation morphism
in the setting of \cite{Nir08-Mod}. The pair is defined as a non vanishing morphism 
$$
\phi:\Gamma\otimes\clE \to \clF
$$
where $\Gamma$ is a finite dimensional $k$-vector space and $\clF$ is a coherent sheaf on $\clX$.
Note that $\phi$ determines a morphism
$$
ev(\phi):\Gamma\otimes\clO_\clX \to \clF\otimes\clE^\vee
$$
which we don't require to be injective on global sections.
Such a definition is reasonable if we think of projective stacks which are banded gerbes.
In that case, there are no morphisms $\Gamma\otimes\clO_\clX\to \clF$ if $\clF$ is not
a pullback from the coarse moduli space. Then in this case by twisting coherent sheaves by the generating  sheaf
we get a richer theory than the theory of the coarse moduli space.
We give a notion of semistability which depends on a Hilbert polynomial of pairs, defined as the sum of
the usual Hilbert polynomial of $\clF\otimes\clE^\vee$ plus a term depending  on a  polynomial $\delta$.
We follow closely two papers dealing with very similar moduli problems:
\cite{HuyLe95} and \cite{Wand10}. In this note we make the exercise of checking
that the proofs extend to our setting, by using results on sheaves over projective stacks
proven in \cite{Nir08-Mod}.
We construct the  stack of  semistable pairs as a global quotient stack and we obtain 
its coarse moduli space with GIT techniques.\\
\subsection*{Conventions}
In this papaer  we work over an algebraically closed field $k$ of characteristic zero. By an {\em algebraic stack} we mean an  algebraic stack over $k$ in the sense of \cite{Art74}. By a {\em Deligne-Mumford stack} we mean an algebraic stack over $k$ in the sense of \cite{DM69}. We assume moreover all stacks and schemes 
unless otherwise stated are are   noetherian  of finite type over $k$.
For sheaves on stacks we refer to \cite{LMBca}.
\subsection*{Notations}
When dealing with sheaves we often adopt the notation in \cite{HuyLe2}.
We denote by $\clX\xrightarrow{\pi} X\to  k$ a projective stack.
We choose a polarization $\clO_X(1)$ on the coarse moduli space. Given a 
sheaf $\clF$ over $\clX$ we often denote by $\clF(m)$ the sheaf $\clF\otimes\pi^*\clO_\clX(m)$.
\section{Introductory material}
\subsection{Recall on projective stacks and on generating sheaves}
\begin{defn}
 A {\em projective stack} is Deligne-Mumford  stack with projective  coarse moduli scheme and  a locally free sheaf  which  is a {\em generating sheaf}  in  the sense of \cite{OlSt03}.
\end{defn}
For the reader's convenience we recall the notion of {\em generating sheaf} following \cite{Nir08-Mod}.
\begin{defn}[{\bf Generating sheaf}]\label{gen-sheaf-defn}
 A locally free sheaf $\clE$ is said to be a {\em generator} for a quasi coherent sheaf $\clF$ is the adjunction morphism (left adjoint to the identity $\pi_*\clF\otimes\clE^\vee \stackrel{\id}{\to}\pi_*\clF\otimes\clE^\vee)$
\begin{eqnarray}
 \theta_\clE(\clF):\pi^*\pi_*\mathcal{H}om_{\clO_\clX}(\clE,\clF)\otimes \clE \to \clF
\end{eqnarray}
is 
surjective. It is a {\em generating} sheaf of $\clX$ if it  is a generator for every quasi coherent sheaf on $\clX$.
\end{defn}

A characterization of generating sheaves can be given by making use of a relative (to the base) ampleness notion for locally free sheaves on stacks introduced in \cite{OlSt03}.
\begin{defn}
 A locally free sheaf  on $\clX$ is {\em $\pi$-ample} if and  only if for every geometric point of $\clX$ the representation on the fiber of the stabilizer group at that point is faithful.
\end{defn}
\begin{defn}
 A locally free sheaf $\clE$ on $\clX $ is   {\em $\pi$-very ample} if for any geometric point of  $\clX$  at that point
the representation  of the stabilizer group on the fiber at that point contains every irreducible representation.
\end{defn}
\begin{prop}[{\cite{Kre06}}, 5.2] Let $\clE$ be a $\pi$-ample sheaf on $\clX$, then there is a  positive integer
$r$ such that the locally free sheav $\bigoplus_{i=0}^r \clE^{\otimes i}$ is $\pi$-very ample.
\end{prop}
\begin{prop}[{\bf \cite{OlSt03}, 5.2}]
 A locally free sheaf on a  Deligne-Mumfors stack $\clX$ is a generating sheaf if and only if it is  $\pi$-very ample.
\end{prop}
We now come to the definition of projective stack.
In \cite{Kre06} it shown that for a proper Deligne-Mumford
stacl over a field the following  characterizations are equivalent.
\begin{thm}[\cite{Kre06}{\bf Corollary 5.4}]
 Let $\clX\to k$ be a proper Deligne-Mumford stack. Then the following are equivalent:
\begin{enumerate}
 \item[1)] the stack $\clX$ has projective coarse moduli space is a quotient stack;
\item[2)] the stack $\clX$ has a projective coarse moduli scheme and there exists a generating sheaf;
\item[3)] the stack $\clX$ has a closed embedding in a smooth proper Deligne-Mumford stack over $k$
and has a projective coarse moduli scheme.
\end{enumerate}
\end{thm}
The third statement is used in \cite{Kre06} as a definition.
\begin{defn}[ {\bf  \cite{Kre06}  Definition 5.5 }]\label{proj-stack-defn}
 A stack $\clX\to k$ is projective if it admits a closed embedding into a smooth Deligne-Mumford stacks proper over $k$
and has projective coarse moduli space.
\end{defn}

\begin{defn}[{\bf Functors $F_\clE$ and $G_\clE$}]\label{functors-defn}
 Let $\clE$ be a locally free sheaf on $\clX$. Let $F_\clE:\mathfrak{QCoh}_{\clX/S}\to \mathfrak{QCoh}_{X/S}$ 
be the functor mapping $\clF\mapsto \pi_*\mathcal{H}om_{\clO_\clX}(\clE,\clF)$ and let $G_\clE:  \mathfrak{QCoh}_{X/S}\to \mathfrak{QCoh}_{\clX/S}$ be a second functor mapping $F\mapsto \pi^*F\otimes\clE$.
\end{defn} 
\begin{rmk} 
 The functor $F_\clE$ is exact because  both  $\otimes \clE^\vee$ and $\pi_*$ are exact functors.
On the other hand $G_\clE$  is not exact unless $\pi^*$ is exact, i.e. $\pi$ is flat.
Examples of stacks with flat map to the coarse moduli scheme are flat gerbes
(e.g. \cite{LMBca} {\bf Definition 3.5}) over schemes or 
stacks root of line bundles (see \cite{AGV06} or \cite{Cadm03}).
\end{rmk}
\begin{rmk}
 The notation $F_\clE$ is the same as in \cite{OlSt03} but $\clG_\clE$ there corresponds to $\clG_E\circ \clF_\clE$ here.
\end{rmk}
\begin{notation}\label{iota-notation}
 We denote by
\begin{eqnarray}
 \iota_\clE(\clF): \clF\otimes\clO_X\rightarrow \clF\otimes\clE nd(\clE)
\end{eqnarray}
the injective morphism mapping a section  to its tensor product with the identity endomorphism of $\clE$.
\end{notation} 
In this paper we use the following notion of slope. Given a coherent sheaf $\clF$ of dimension $d$
\begin{eqnarray}
 \hat{\mu}_\clE(\clF)=\frac{\alpha_{d-1}(\clF\otimes\clE)}{\alpha_d(\clF\otimes\clE^\vee)}.
\end{eqnarray}
We often denote the multiplicity $\alpha_d(\clF\otimes\clE^\vee)$ by $r_{\clE,\clF}$.

\section{Setting up  the moduli problem}

In the following $\pi:\clX\to X \to k$ will be a smooth projective stack  with coarse moduli scheme $X$ over an algebraically closed field  $k$.  
We will fix a polarization $(\clO_X(1),\clE)$ and a rational polynomial $\delta$ such that $\delta(m)\geq 0$ for $m>>0$.
\begin{defn}\label{pairs-defn}
 A {\em  pair}  $(\clF,\phi)$ is a non-trivial morphism
\begin{eqnarray}
 \phi:\Gamma\otimes \clE\to \clF,
\end{eqnarray}
where $\Gamma$ is a finite dimensional $k$-vector space,  $\clF$ is a coherent sheaf of dimension $d$, $d\in\bbN$, $d\leq \mbox{dim}\ \clX$, 
and $\clE$ is the fixed generating sheaf.  A morphism between two pairs $(\clF,\phi)$, $(\clF',\phi')$ is 
a commutative diagram
\begin{eqnarray}\label{morphs-diag-def}
 \xymatrix{
\Gamma\otimes\clE\ar[r]^{\phi}\ar[d]_{\lambda} & \clF\ar[d]^{\alpha}\\
\Gamma\otimes \clE\ar[r]_{\phi'} & \clF'
}
\end{eqnarray}
where $\lambda\in \bbC^*$ and $\alpha$ is a morphism of coherent sheaves.
\end{defn}

Note that we can relate the notion of pairs to a {\em stacky version} of coherent systems of \cite{LePot93}.
For the reader's convenience we recall the definition of coherents systems on schemes given by LePotier.
\begin{defn}{\bf \cite{LePot93} Def. 4.1}
 Let $X$ be a smooth projective variety of dimension $n$. A coherent system of dimension $d$ is a pair
$(\Gamma,F)$, where $F$ is a coherent sheaf of dimension $d$ over $X$ and $\Gamma\subseteq H^0(F)$ is a vector subspace.
\end{defn}
We extend this notion to projective Deligne-Mumford stacks.
\begin{defn}
 Let $\clX$ be  a projective Deligne-Mumford stack over $k$. A {\em twisted coherent system} on $\clX$ is a pair
$(\Gamma,\ \clF)$ where $\clF$ is a coherent sheaf over $\clX$  and $\Gamma\subseteq H^0(\clF\otimes\clE^{\vee})$ 
is a vector subspace.
\end{defn}

A pair $(\clF,\phi)$  determines a subspace   of  $H^0(\clX,\clF\otimes\clE^\vee)$ given
by the image of $\Gamma$ along $H^0(ev(\phi))$  where 
\begin{eqnarray}
 \xymatrix{
ev(\phi):\Gamma\otimes\clO_\clX\ar@{^{(}->}[r]^-{\iota_\clE(\Gamma)} & (\Gamma\otimes\clE)\otimes\clE^\vee\ar[r]^-{\phi\otimes\clE^\vee} & \clF\otimes\clE^\vee
}
\end{eqnarray}
is obtained from $\phi$ by applying the functor $-\otimes\clE^\vee$ and by composing by the inclusion morphism $\iota_\clE(\Gamma)$.
Note that $H^0(ev(\phi))$ is  not necessarily injective on global sections. Hence
the twisted coherent system determined as above is $(W,\clF)$ with $\mbox{dim} W < \mbox{dim} \Gamma$.
Conversely, let  us consider a twisted pair $(\Gamma,\clF)$. Let 
\begin{eqnarray}\label{descr-1}
 ev: \Gamma\otimes\clO_\clX\to \clF\otimes\clE^\vee.
\end{eqnarray}
be the corresponding evaluation morphism.
It  is  possible to associate to (\ref{descr-1})
  the  pair $(\clF,\phi)$
\begin{eqnarray}\label{induced-ev-morphi-eqn}
 \xymatrix{
\phi(ev): \Gamma\otimes\clE\ar[r]^-{ev\otimes\clE} & (\clF\otimes \clE^\vee)\otimes \clE\ar[r]^-{Tr} & \clF
}
\end{eqnarray}
obtained by applying the functor  $-\otimes\clE^\vee$ and by composing with $Tr:\clE nd(\clE)\to\clO_\clX$. 
It is not hard to see that $ev(\phi(ev))=ev$ and that $\phi(ev(\phi))=\phi$.

We also give the definition of family  of pairs.
Let $S$ be a scheme of finite type over $k$. Let $\pi_\clX:\clX\times S\to \clX$ and
$\pi_S:\clX\times S\to S$ be the natural projections. 
\begin{defn}\label{family-of-pair-defn}
 A {\em pair parametrized  by $S$}  is a $S$-flat coherent 
sheaf $\clF$ over $\clX\times S$ and a homomorphism
$$
\phi_S:\pi_\clX^* \Gamma\otimes\clE\to \clF 
$$
such that for any closed point $s$ of $S$ 
$$
\phi_S(s):\pi_\clX^* \Gamma\otimes\clE(s)\to \clF(s) 
$$
is a pair.
\end{defn}

\subsection{(Semi)stability}
We define a parameter-dependent  Hilbert polynomial for a  stable pair in the following way.
\begin{defn}\label{hilb-polyn-defn}
 The {\em Hilbert polynomial} of a pair $(\clF,\phi)$ is 
\begin{eqnarray}
 P_\clE{(\clF,\phi)}:=P(F_\clE(\clF)) +\epsilon(\phi)\delta
\end{eqnarray}
where $\epsilon(\phi)=1$ if $\phi\neq 0$ and $0$ otherwise. \\
The {\em reduced Hilbert polynomial} is 
\begin{eqnarray}
p_\clE(\clF,\phi):= \frac{P_\clE{(\clF,\phi)}}{r_{F_\clE(\clF)}}.
\end{eqnarray}
\end{defn}
\begin{rmk}
Note that in the above  definition following \cite{Nir08-Mod} we do not use the Hilbert polynomial
$P(\clF)$ of the sheaf on the stack, rather the Hilbert polynomial  $P(F_\clE(\clF))=P(\clF\otimes \clE^\vee)$.
The reason is particularly evident if we consider sheaves on gerbes. In this case the non twisted Hilbert polynomial
of any sheaf which is not a pull-back from the coarse moduli space vanishes. For more details see \cite{Nir08-Mod}.
\end{rmk}
We will define (semi)stability by using the Hilbert polynomial introduced above. We need some more preliminary remarks and notations.
\begin{defn}\label{induced-hom-defn}
 Let $(\clF,\phi)$ be a stacky pair. Any subsheaf $\clF'\subset\clF$ defines a {\em induced homomorphism}
$\phi':\Gamma\otimes\clE\to \clF'$ which is   equal to $\phi$ if $\mbox{Im}\ \phi\subseteq \clF'$ and zero otherwise.
The corresponding quotient $\clF''=\clF/\clF'$ also inherits  an {\em induced homomorphism}
$\phi'':\Gamma\otimes \clE\to\clF''$ which is defined as the composition of $\phi$ with the quotient map. Note that it is the zero morphism if and only if $\mbox{Im}\ \phi\subseteq \clF'$.
\end{defn}
\begin{rmk}
 The Hilbert polynomial of stacky pairs is additive on short exact sequences.
\end{rmk}

\begin{defn}
 A stacky pair is {\em  $\delta$ (semi)stable} if for any saturated submodule $\clF'\subset \clF$
\begin{eqnarray}\label{delta-ss-defn-eqn}
 P_\clE(\clF',\phi') (\leq) r_{\clE,\clF'}\  p(F_\clE(\clF),\phi)
\end{eqnarray}
\end{defn}

\begin{defn}\label{ss-pair-family-defn}
 A {\em $\delta$ semi stable   pair} parametrized by $S$ is a stacky pair over $S$ such that for every closed  point of $S$ the pair
$(\clF(s), \phi|_{\pi_\clX^*\Gamma\otimes\clE_(s)})$ is a $\delta$ semistable  pair.
\end{defn}

\subsection{Properties of $\delta$ semi stable  pairs.}
We note that $\delta$ semi stablity implies purity of the underlying sheaf  of the pair. 
\begin{prop}
 Let $(\clF,\phi)$ be  a $\delta$ semi stable pair. Then  $\clF$ is pure.
\end{prop}
\begin{pf}
 Let us assume that $\clF$ is not pure. Let $\mathcal{T}=T_{d-1}(\clF)$ be  the element of the torsion filtration of maximal dimension. Then
\begin{eqnarray}
 P_\clE(\mathcal{T})+\epsilon(\mathcal{T})\leq \frac{r_{F_\clE(\mathcal{T})}}{r_{F_\clE(\clF)}}(p+\delta)=0
\end{eqnarray}
where the equality on the r.h.s. holds because $\mathcal{T}$ is a sheaf at most  of dimension $d-1$.
It follows that $\clT=0$.
\end{pf}
For $\delta$ semistability we have  characterizations and properties analogous to usual 
Gieseker semistability for sheaves on schemes. We list some of them
\begin{prop}
Let $(\clF,\phi)$  be a stacky pair. Then the following conditions are equivalent:
\begin{enumerate}
\item[i)]for all proper  subsheaves $\clF'\subseteq \clF$ $$P_\clE(\clF',\phi') (\leq) r_{F_\clE(\clF')}p_\clE(\clF,\phi);$$
 \item[ii)] $(\clF,\phi)$ is (semi)stable;
\item[iii)] for all proper quotient sheaves $\clF\to \clF''$ with $\alpha_d(F_\clE(\clF))>0$  $$P_\clE(\clF'',\phi'')(\geq) r_{F_\clE(\clF'')} p_\clE(\clF,\phi);$$
\item[iv)] for all proper purely $d$-dimensional quotient sheaves $\clF\to \clF''$ with $\alpha_d(F_\clE(\clF))>0$  $$P_\clE(\clF'',\phi'') (\geq) F_\clE(\clF'') p_\clE(\clF,\phi);$$
\end{enumerate}
\end{prop}
\begin{pf}
 The proof is very similar to \cite{HuyLe2} {\em Prop. 1.2.6}. We use additivity of the ranks and of the modified Hilbert polynomials on short exact sequences. Note that 
  if  inequality
(\ref{delta-ss-defn-eqn}) holds for saturated subsheaves, it also holds for arbitrary subsheaves.
Indeed, let $\clF'\subseteq \clF$ be a not necessarily saturated subsheaf. Then if $\mbox{Im}\ phi$ is contained in 
$\clF'$, then $\mbox{Im}\ \phi\subseteq \clF'^s$, where $\clF'^s$ is the saturation of $\clF'$ in $\clF$. Moreover
$P_\clE(\clF')\leq P_{\clE}(\clF'^s)$.
\end{pf}
\begin{lem}\label{ss-hom-lem}
 Let $(\clF,\phi)$, $(\clG,\psi)$  be  two $\delta$ semistable pairs such that 
$p_\clE((\clF,\phi))> p_\clE((\clG,\psi))$. Then $\mbox{Hom}(\clF,\phi),(\clG,\psi))=0$.
\end{lem}
\begin{pf}
 Let us assume there is a non zero morphism  $(\alpha,\lambda)$. Let $\clH=\mbox{Im}\ \alpha$.
 By semi stability we get
\begin{eqnarray}\label{chain-ineq-eqn}
 p_\clE(\clF,\phi)\leq p_\clE(\clH,\phi_\clH)= p_\clE(\clH,\psi_\clH)\leq p_\clE(\clG,\psi),
\end{eqnarray}
where $\phi_\clH$ and $\psi_\clH$ are the induced homomorphisms. Inequality  (\ref{chain-ineq-eqn})
contradicts the assumption.
\end{pf}

\begin{lem}
 Let $\alpha: (\clF,\phi)\to (\clG,\psi)$ be a homomorphism between $\delta$ stable pairs of the same reduced Hilbert polynomial.
Then $\alpha$ is 0 or an isomorphism.
\end{lem}
\begin{pf}
 Analogous to {\em Lemma 1.6} in \cite{Wand10}. Cfr. also \cite{HuyLe2} {\em Proposition 1.2.7}.
\end{pf}
\begin{cor}
 Let $(\clF,\delta)$ be a semi stable pair. Then $\mathcal{E}nd((\clF,\delta))$ is a finite dimensional division algebra. Since we work over an algebraically closed field $k$ $\mathcal{E}nd((\clF,\delta))\simeq k$.
\end{cor}
\begin{pf}
 Same as \cite{HuyLe2} {\em Cor. 1.2.8}.
\end{pf}
\subsection{Harder-Nahrasiman and Jordan-H\"older filtration}
\begin{prop}
 Let $(\clF,\phi)$ be a  pair such that $\clF$ is pure. Then it admits a unique  Harder-Nahrasiman filtration 
\begin{eqnarray}
 0\subset HN_0(\clF,\phi)\subset....\subset HN_{l-1}(\clF,\phi)\subset HN_{l}(\clF,\phi)=(\clF,\phi)
\end{eqnarray}
such that each $gr_i^{HN}(\clF,\phi)= HN_i(\clF,\phi)/HN_{i-1}(\clF,\phi)$ is $\delta$ semistable
and if $p_i:=p_\clE(gr_i^{HN}(\clF,\phi))$ then
\begin{eqnarray}
p_{max}(\clF,\phi)=p_1 > p_2 > ...>p_l=p_{min}(\clF,\phi)
\end{eqnarray}
\end{prop}
\begin{pf}
 The proof proceeds as in \cite{HuyLe2} Theorem 1.3.4. 
Indeed, it is possible to find a subsheaf $\clF_0\subseteq \clF$ such that
it is not contained in any subsheaf $\clF'$  of $\clF$ with 
$p_\clE(\clF_0,\phi_0) < p_\clE(\clF',\phi')$, where $\phi_0$ and $\phi'$ are the induced homomorphisms.
This implies the existence part. Uniqueness is proven by using Lemma \ref{ss-hom-lem}.
\end{pf}

\begin{prop}
 Let $(\clF,\phi)$ be a $\delta$ semistable pair with reduced Hilbert polynomial $p$. Then there is a  Jordan-Holder filtration
\begin{eqnarray}
0= JH_0(\clF,\phi)\subset JH_1(\clF,\phi)\subset ...\subset JH_l(\clF,\phi)=(\clF,\phi)
\end{eqnarray}
such that each $gr_i^{JH}(\clF,\phi)= JH_i(\clF,\phi)/JH_{i-1}(\clF,\phi)$ is $\delta$ stable with reduced Hilbert polynomial $p$.
The graded object $gr^{JH}(\clF,\phi)=\oplus_i\  gr_i^{JH}(\clF,\phi)$ is independent of the choice of the filtration. 
Note that  it inherits an induced homomorphism $gr^{JH}(\phi):\Gamma\otimes\clE\to gr^{JH}(\clF,\phi)$.
\end{prop}
\begin{pf}
 The proof is the same as \cite{HuyLe2} {1.5}. The same arguments hols because of additivity  the modified Hilbert polynomial on short exact sequences.
\end{pf}
\begin{rmk}
 It is not hard to see that $gr^{JH}(\phi)$ is non trivial if $\phi$ is not, and its image is cointained in only one summand of $gr^{JH}(\clF,\phi)$.
\end{rmk}
\begin{defn}\label{S-equiv-defn}
 Two $\delta$ semistable pairs are said to be {\em $S$-equivalent} if their Jordan-H\"older graded objects are isomorphic.
\end{defn}

\begin{rmk}
 From now on we  will always assume that $\delta$ is strictly positive and $\phi$ is non vanishing
otherwise $\delta$ (semi)stability reduces to usual Gieseker (semi)stability for sheaves on stacks.
\end{rmk}

We introduce  a symbol which is convenient to restate the (semi)stability
 condition when assuming that the homomorphism of a pair is non vanishing.
\begin{defn}\label{epsilon-defn}
 Let $(\clF,\phi)$ be a stacky pair. For any exact sequence
\begin{eqnarray}
 0\to \clF'\to \clF\to \clF''\to 0
\end{eqnarray}
let
\begin{eqnarray}
 \epsilon(\clF') := \left\{ \begin{array}{rl} 1 & \quad \mbox{if Im}\ \phi\subseteq \clF' \\
                                          0 & \quad \mbox{otherwise}
                         \end{array}\right.\nonumber
\end{eqnarray}
and
\begin{eqnarray}
\epsilon(\clF''):= 1 -\epsilon(\clF')
\end{eqnarray}
\end{defn}
With the above definition we can restate the $\delta$ (semi)stability condition for stacky pairs $(\clF,\phi)$. Indeed $(\clF,\phi)$
is semi stable if and only if for every saturated subsheaf $\clF'$
\begin{eqnarray}\label{nice-notation-sst-restatement}
 P_\clE(\clF') + \epsilon(\clF')\delta (\leq) \frac{r_{F_\clE(\clF')}}{r_{F_\clE(\clF)}} (P_\clE(\clF) + \delta)
\end{eqnarray}

\section{Boundedness}\label{bounded-prop}
In this section we prove boundedness of the family of  $\delta$ semistable  pairs.
Since we want to use a GIT construction similar to  \cite{HuyLe95} and \cite{Wand10} we take $\mbox{deg}\ \delta < \mbox{dim} \clX$.
\begin{prop}
 Let $P$ be  a fixed  polynomial of degree $d< \mbox{dim}\ \clX$. Then the family of $\delta$ semistable  pairs with 
 Hilbert polynomial $P$ is bounded.
\end{prop}
\begin{pf}
Let $\mathfrak{F}$ be a family of  coherent sheaves over $\clX$.
According to \cite{Nir08-Mod} Corollary 4.17 $\mathfrak{F}$ is bounded if and only if $F_\clE(\mathfrak{F})$ is bounded over $X$.
 We use \cite{Simp-rep-I} Theorem 1.1,  according to which a family $\mathfrak{F}$  of sheaves over a  projective scheme  $X$   with fixed Hilbert polynomial  is bounded if and only there exists a constant $C$ such that for any $F\in \mathfrak{F}$  $\mu_{max}(F)\leq C$.
Let $\phi: \clE\otimes \Gamma\to \clF$ be a  $\delta$ semistable pair. Let $\mbox{Supp}(\clF)=\clY$. 
Let us assume first that $\mbox{Im}\phi\nsubseteq HN_{l-1}(\clF)$, where $HN_{l-1}(\clF)$ is the 
maximal proper subsheaf in the Harder-Nahrasiman filtration.
 Then the composition
\begin{eqnarray}
\Gamma\otimes \clE \otimes\clO_\clY\to  \clF\to gr^{HN}_l(\clF) 
\end{eqnarray}
is a  non zero morphism between sheaves of pure dimension $d$. This implies that
\begin{eqnarray}\label{mu-min-eqn}
 \hat{\mu}_{\clE,min}(\clE \otimes\clO_\clY)\leq \hat{\mu}_{\clE, min}(\clF).
\end{eqnarray} 
We note first that  $\hat{\mu}_{\clE,min}(\clE \otimes\clO_\clY)\geq \hat{\mu}_{min}(\pi_*\clE nd(\clE)\otimes\clO_Y)$
where $Y=\pi(\clY)$. The reason is that not all the quotient sheaves of  $F_\clE(\clF)$  are obtained as images by the functor $F_\clE$
of quotient sheaves of $\clF$ (cfr.\cite{Nir08-Mod} {\bf Remark 3.15})
We want to find a lower bound for $\mu_{min}(\clO_Y)$. This is provided by a result proven in
 \cite{LePot93}. 
\begin{cor}[{\bf\cite{LePot93} Corollary 2.13} ]
 Let $X$ be a projective scheme over $k$. Let $S$ be a subscheme of pure dimension $d$ and of degree $k$. Then $\mu_{min}(\clO_S)$ is bounded from below by a constant which only depends on
$d$, $k$ and $X$.
\end{cor}
We observe that $Y=\pi\ \mbox{Supp}\ \clF=\mbox{Supp}\ F_\clE(\clF)$ is a  purely $d$-dimensional subscheme
of degree $\leq r_{\clE,\clF}^2$.
Then 
\begin{eqnarray}\label{mu-min-bound-eq }
\hat{\mu}_{\clE,min}(\clF)\geq \hat{\mu}_{min}(\clO_Y)+ \hat{\mu}_{min}(\pi_*\clE nd(\clE))\geq A + \hat{\mu}_{min}(\pi_*\clE nd(\clE)):=B
 \end{eqnarray}
for some constant $A$ which only depends on $X$ and on the fixed polinomial $P$.
 By the barycenter formula for the slope this implies that
\begin{eqnarray}
 \hat{\mu}_{\clE,max}(\clF)\leq \mbox{max}\{r_{\clE,\clF}\ \hat{\mu}_\clE(\clF) - (r_{\clE,\clF}-1)\ B,  \hat{\mu}_\clE(\clF) \}.
\end{eqnarray}
Boundedness of $\mu_{max}(F_\clE(\clF))$ follows from Lemma \ref{bound-on-stack-implies-cms}.
Let us consider now the case where $\mbox{Im}\phi\subseteq HN_{l-1}(\clF)$.
Then by $\delta$ semistability
$$
p_\clE(HN_{l-1}(\clF))\leq p_\clE(\clF),
$$
which in turn implies that
$$
p_\clE(\clF)\leq p_\clE(gr_l(\clF))
$$
and by the barycenter formula that
$$
\hat{\mu}_{\clE,max}(\clF)\leq \hat{\mu}_\clE(\clF).
$$
Summing up we get
\begin{eqnarray}\label{bound-on-mu-max-final}
 \hat{\mu}_{max}(F_\clE(\clF))\leq \mbox{max}\{\hat{\mu}_\clE(\clF) ,r_{\clE,\clF}\mu_\clE(\clF) - (r_{\clE,\clF}-1) B \} + \tilde{m}\  \mbox{deg} \clO_X(1),
\end{eqnarray}
where the inequality is a consequence of the above estiamtes and of {\em Lemma \ref{bound-on-stack-implies-cms}}. 
\end{pf} 

\begin{lem}\label{bound-on-stack-implies-cms}
 Let $\clF$ be a coherent sheaf on $\clX$ of pure dimension $d$. Then
if $\mu_{max,\clE}(\clF)$ ($\mu_{min,\clE}(\clF)$) is bounded from above (below), then also $\mu_{max}(F_\clE(\clF))$ ($\mu_{max}(F_\clE(\clF))$) is bounded from above (below). 
\end{lem}
\begin{pf}
The proof is similar to \cite{Nir08-Mod} {\bf Proposition 4.24}.
 Let  $\overline{F}\subset F_\clE(\clF)$ be the maximal destabilizing subsheaf.
Let us consider the morphism
\begin{eqnarray}
 \pi^*\overline{F}\otimes \clE \longrightarrow \pi^*(\pi_*\clF\otimes\clE^\vee)\otimes \clE \xrightarrow{\theta_\clE(\clF)} \clF.
\end{eqnarray}
The right arrow is surjective by definition of generating sheaf. Let $\overline{\clF}$ be the subsheaf corresponding to the image 
of the composition.  By applying the functor $F_\clE$ we get  the  surjective morphism
\begin{eqnarray}
 \overline{F}\otimes \pi_*\mathcal{E}nd(\clE)\to F_\clE(\overline{\clF})
\end{eqnarray}
Let $\tilde{m}>>1$ be an integer number   such that $\pi_*\mathcal{E}nd(\clE)$ is generated by global sections.
Let $N=h^0(\pi_*\mathcal{E}nd(\clE)(\tilde{m}))$. Then there is a surjective morphism
$\overline{F}\otimes\clO_X^{\otimes N}(-\tilde{m})\to F_\clE(\overline{\clF})$. Note that since $\overline{F}$ is semistable, so is
 $\overline{F}(-\tilde{m})$. Moreover, for any $k\in\bbN$, $\overline{F}(-\tilde{m})^{\oplus k}$  is also semistable. 
By composition we get the surjective morphism
\begin{eqnarray}
\overline{F}(-\tilde{m})^{\oplus N} \to \overline{F}\otimes \pi_*\mathcal{E}nd(\clE)\to F_\clE(\overline{\clF}).
\end{eqnarray}
Then by semistability of $\overline{F}(-\tilde{m})^{\oplus N}$ 
\begin{eqnarray}
\hat{\mu}_{max}(F_\clE(\clF))\leq \hat{\mu}_\clE(\overline{\clF}) + \tilde{m}\ \mbox{deg}\ \clO_X(1)
\leq \hat{\mu}_{\clE,max}(\clF) + \tilde{m}\ \mbox{deg}\ \clO_X(1).
\end{eqnarray}
Hence  if  $\hat{\mu}_{\clE,max}(\clF)$ is bounded from above, $\hat{\mu}_{max}(F_\clE(\clF))$ is also bounded from above.
Boundedness from below of $\hat{\mu}_{max}(F_\clE(\clF))$ also follows.
\end{pf}

\subsection{Rephrasing semistability in terms of number of global sections}
We apply here a result due  Le Potier and Simpson (cfr. e.g. \cite{HuyLe2} Corollary 3.3.1 and 3.3.8) in order
to get a bound  on the number of global sections of  $F_\clE(\clF)$, where $\clF$ is a coherent
sheaf on $\clX$ of pure dimension $d$.
 We state it for sheaves on $X$ obtained by applying the functor $F_\clE$ to some sheaf on $\clX$.
\begin{cor}\label{glob-sections-bound-char0-coro}
 Let $\clF$ be a $d$-dimensional coherent sheaf over $\clX$. Let $r=r_{\clE,\clF}$ be the 
multiplicity of $F_\clE(\clF)$.  Let $C:=r(r+d)/2$. Then
\begin{displaymath}
 h^0(\clF\otimes\clE^\vee)\leq \frac{r-1}{r}\cdot\frac{1}{d!}[\hat{\mu}_{max}(F_\clE(\clF))+C-1+m]^d_+ +\frac{1}{r}\cdot
\frac{1}{d!}[\hat{\mu}(F_\clE(\clF))+C-1+m]^d_+,
\end{displaymath}
whre $[x]_+=\mbox{max}\ \{0,x\}$. 
\end{cor}

We will need the above estimate  in order to
give the following characterization of semistability.

\begin{prop}\label{num-cond-d-ss-prop}
 For $m>>0$ for any pure pair $(\clF,\phi)$ the following properties are equivalent  
\begin{enumerate}
 \item $(\clF,\phi)$ is $\delta$ (semi)-stable,
\item $P(m)\leq h^0(\clF\otimes\clE^\vee(m))$ and for any subsheaf $\clF'\subseteq \clF$ with $0<r_{F_\clE(\clF')}<r_{F_\clE(\clF)}$
$$
h^0(\clF'\otimes\clE^\vee(m))+\epsilon(\clF')\delta(m) (\leq) < \frac{r_{F_\clE(\clF')}}{r_{F_\clE(\clF)}}(P(m)+\delta(m))
$$ 
\item for any quotient $\clF\to \clF''$ with $0<r_{F_\clE(\clF'')}<r_{F_\clE(\clF)}$
$$
\frac{r_{F_\clE(\clF'')}}{r_{F_\clE(\clF)}}(P(m)+\delta(m))(\leq) < h^0(\clF''\otimes\clE^\vee(m))+\epsilon(\clF'')\delta(m)
$$
\end{enumerate}
\end{prop}
\begin{pf}
 We prove $(1)\Rightarrow (2)$. The family of sheaves
underlying the family of $\delta$ semistable pairs on $\clX$ is bounded.
This is equivalent to the family of sheaves on $X$ obtained by applying the functor $F_\clE$ being bounded
(cfr. \cite{Nir08-Mod} {\bf Corollary 4.17}).
Hence there exists $m$ such that
for any $\clF$ in a $\delta$ semistable pair  $P_\clE(m)=h^0(F_\clE(\clF)(m))$. Let $\clF'\subseteq\clF$
be an arbitrary subsheaf of $\clF$. Let us assume that inequality (\ref{bound-on-mu-max-final})
gives $\hat{\mu}_{max}(F_\clE(\clF))\leq \hat{\mu}_{max}(F_\clE(\clF)) + \tilde{m}\mbox{deg}\ \clO_X(1)$.
 We distinguish two cases:\\
\begin{enumerate}\label{cases-distinction}
 \item[A)] $\hat{\mu}(F_\clE(\clF'))\geq \hat{\mu}_\clE(\clF) - (r-1)\tilde{m}\mbox{deg}\ \clO_X(1) - C\cdot r - \delta_1 r$;
\item[B)]  $\hat{\mu}(F_\clE(\clF'))\leq \hat{\mu}_\clE(\clF) -(r-1)\tilde{m}\mbox{deg}\ \clO_X(1) - C\cdot r - \delta_1 r$;
\end{enumerate}
where $C=r(r+d)/2$ and $r=r_{\clE,\clF}$. If $\clF'$ is of type $A$, then $\hat{\mu}(F_\clE(\clF'))$ is bounded from below.  We observe that we can assume that 
$\clF'$ is saturated (which implies that also $F_\clE(\clF')$ is saturated, because the functor 
$F_\clE$ maps torsion filtrations to torsion filtrations, cfr. \cite{Nir08-Mod} {\bf Corollary 3.17}). Indeed for any sheaf $\clH$ on $\clX$
 $P_\clE(\clH)\leq P(\clH^{s})$,
where $\clH^{s}$ is the saturation of $\clH$. Then the family of sheaves of type $A$ is bounded by Grothendieck Lemma for  stacks
(see \cite{Nir08-Mod} {\bf Lemma 4.13}). As a consequence the number of Hilbert polynomials of the family is finite and 
there exists an integer number $m_0$ such that for  any $m\geq m_0$ and for any subsheaf $\clF'$ of type $A$
$P(\clF'\otimes\clE^\vee)=h^0(\clF'\otimes\clE^\vee)$ and
\begin{eqnarray}
 P(\clF'(m)\otimes\clE(m)) + \epsilon(\clF')\delta(m) & (\leq) & \frac{r_{F_\clE(\clF')}}{r_{F_\clE(\clF)}}[p(m)+\delta(m)]\Leftrightarrow\nonumber\\  P(\clF'\otimes\clE^\vee) + \epsilon(\clF')\delta & (\leq) & \frac{r_{F_\clE(\clF')}}{r_{F_\clE(\clF)}}[p+\delta].
\end{eqnarray}
Let us consider    now sheaves of type $B$. We get
\begin{displaymath}
 h^0(\clF'\otimes\clE^\vee)\leq \frac{r'-1}{r'}\cdot\frac{1}{d!}[\hat{\mu}_{max}(F_\clE(\clF'))+C'-1+m]^d_+ +\frac{1}{r'}\cdot
\frac{1}{d!}[\hat{\mu}(F_\clE(\clF'))+C'-1+m]^d_+,
\end{displaymath}
where $C'=r'(r'+d)/2$ and $r'=r_{\clE,\clF'}$. This in turn implies
\begin{eqnarray}
 \frac{h^0(\clF'\otimes\clE^\vee)}{r_{\clE,\clF'}} & \leq  &  \frac{r-1}{r}\cdot\frac{1}{d!}[\hat{\mu}_\clE(\clF) +\tilde{m}\mbox{deg}\ \clO_X(1) +C-1+m]^d_+ \nonumber \\
& & +\frac{1}{r}\cdot
\frac{1}{d!}[\hat{\mu}(F_\clE(\clF))+ (1-r)C -(r-1)\tilde{m}\mbox{deg}\ \clO_X(1) -1 -\delta_1r + m]^d_+\nonumber\\
& \leq & \frac{m^d}{d!} + \frac{m^{d-1}}{(d-1)!}(\hat{\mu}_\clE(\clF) -1 -\delta_1)+\ .....
\end{eqnarray}
where $\delta_1$ is the degree $d-1$ coefficient of $\delta$ and 
 $....$ stay for lower degree polynomials.
We can conclude that
\begin{eqnarray}
 \frac{1}{r_{\clE,\clF'}}(h^0(\clF'\otimes\clE^\vee(m))+\epsilon(\clF')\delta(m))
< \frac{P(m)}{r} <  \frac{P(m)+\delta(m)}{r}
\end{eqnarray}

If inequality (\ref{bound-on-mu-max-final}) gives a different upper bound 
it is possible to apply the same arguments used above, except that suitably modified bounds have to be chosen
to define sheaves of type $A$ and $B$ as on page
 \pageref{cases-distinction}.\\
$(2)\Rightarrow (3)$
 Let $\clF''$ be any quotient of $\clF$ with $0<r_{F_\clE(\clF'')}<r_{F_\clE(\clF)}$.
Let $\clF'\subseteq\clF$ denote the corresponding kernel.
Then
\begin{eqnarray}
 h^0(\clF''\otimes\clE^\vee(m))+\epsilon(\clF'')\delta(m) & (\geq) & h^0(\clF\otimes\clE^\vee(m))-h^0(\clF''\otimes\clE^\vee(m))+\delta(m)-\epsilon(\clF'')\delta(m)\nonumber\\
&(\geq)& \frac{1}{r}(rP_\clE(m)-r'P_\clE(m)+r\delta(m)-r'\delta(m))\nonumber\\
& = & \frac{r''}{r'}(P_\clE(m)+\delta(m))\nonumber
\end{eqnarray}
whhere $r=r_{F_\clE(\clF)}$, $r'=r_{F_\clE(\clF')}$, $r''=r_{F_\clE(\clF'')}$.\\
$(3)\Rightarrow (1)$.  
Let us show first that the underlying sheaves of $\delta$ semistabe pairs satisfying $(3)$ form a bounded family.
Let $(\clF,\phi)$ be  a such a   pair. 
 Let ${\clF}_{min}$  the minimal destabilizing quotient   of $\clF$. Then by hypothesis
\begin{eqnarray}
 \frac{P_\clE(m)+\delta(m)}{r_{\clE,\clF}} - \frac{\epsilon(\clF_{min})\delta(m)}{r_{\clE,\clF_{min}}}\leq
  \frac{h^0(\clF_{min}\otimes \clE^\vee)}{r_{\clE,\clF_{min}}} \nonumber\\
\leq \frac{1}{d!}[\hat{\mu}_\clE(\clF_{min})+\tilde{m}\mbox{deg}\ \clO_X(1) +C -1 +m]_+^d \nonumber
\end{eqnarray}
where the second inequality can be deduced from {\em Corollary} \ref{glob-sections-bound-char0-coro}.
It follows that $\hat{\mu}_{\clE,min}(\clF)$ is bounded from below.  By Lemma \ref{bound-on-stack-implies-cms}, $\hat{\mu}_{min}(F_\clE(\clF))$ is also bounded from below or equivalently $\hat{\mu}_{max}(F_\clE(\clF))$ is bounded from above, 
which implies boundedness for sheaves satisfying $(3)$.
 Let us consider now an arbitrary quotient
$\clF''$ of $\clF$. Then either $\hat{\mu}_{\clE}(\clF'')>\hat{\mu}_{\clE}(\clF)+\delta_1/r$ or $\hat{\mu}_{\clE}(\clF'')\leq\hat{\mu}_{\clE}(\clF)+\delta_1/r$. In the first case  a strict inequality for Hilbert polynomials is satified, implying stability.
 In the second case  $\hat{\mu}_{\clE}(\clF'')$ is bounded from above.  By Grothendieck Lemma for stacks
we get that the family of pure dimensional quotients satisfying the second inequality is  bounded. Hence there exists  $m$ large enough
such that for any sheaf $\clF$ satisfying $(3)$ 
$h^0(\clF''\otimes\clE^\vee(m))=P_\clE(\clF''(m))$. Therefore
\begin{eqnarray}
 P(\clF''\otimes\clE^\vee(m))+\epsilon(\clF'')\delta(m)(\geq)\frac{r''}{r}[P_\clE(m)+\delta(m)]\nonumber\\\Leftrightarrow P(\clF''\otimes\clE^\vee)+\epsilon(\clF'')\delta (\geq) \frac{r''}{r}[P_\clE+\delta],
\end{eqnarray}
where $r=r_{\clE,\clF}$ and $r''=r_{\clE,\clF''}$.
Eventually we remark  that  from the proofs of $(2)$ and $(3)$  equality holds if and only if the subsheaf or
the quotient sheaf is destabilizing.
\end{pf}  
We recall a technical result originally due to  Le Potier which will be used in Theorem \ref{GIT-implies-dss-thm}.
The proof can be founf in \cite{Nir08-Mod} Lemma 6.10  and generalizes \cite{HuyLe2} Proposition 4.4.2.
\begin{lem}\label{deform-sheaf-lem}
 Let $\clF$ be a coherent sheaf over $\clX$ that can be deformed to a sheaf of the same dimension $d$. Then there is a
pure $d$-dimensional sheaf $\clG$ on $\clX$ with a map $\clF\to \clG$ such that the kernel is
$T_{d-1}(\clF)$ and $P_\clE(\clF)=P_\clE(\clG)$.
\end{lem}

\section{The parameter space}
By {\bf Proposition \ref{bounded-prop}} we know that  the family of sheaves obtained by 
applying the functor $F_\clE$ to the underlying family of sheaves  of $\delta$ semistable pairs with
fixed Hilbert polynomial is bounded.
Therefore there exists an integer number $m_0$ such that for any $m\geq m_0$ and for any$\delta$ semistable pair $(\clF,\phi)$,
$F_\clE(\clF)(m)$ is generated by global sections. Let $V$ be a $k$ vector space of dimension equal to $P_\clE(\clF)(m)$.
There is  a surjective morphism
\begin{eqnarray}\label{q-constr-eqn}
\xymatrix{
q: V\otimes\clE(-m)\ar[r] & \pi^*\pi_*(\clF\otimes\clE^\vee)\otimes\clE \ar@{->>}[r] &  \clF 
}
\end{eqnarray}
obtained by applying the functor $G_\clE$ to 
$$
V\isomto H^0(F_\clE(\clF(m))\twoheadrightarrow F_\clE(\clF(m)
$$
and composing with $\theta_\clE(\clF): G_\clE(F_\clE(\clF))\to \clF$.
The morphism $q$ corresponds to a closed point of $\tilde{Q}:=Quot(V\otimes\clE(-m), P_\clE(\clF))$
(existence of Quot schemes on Deligne-Mumford stacks follows from \cite{OlSt03}).
Up to changing $m$ we can assume that also $F_\clE(\clE)(m)=\pi_*\mathcal{E}nd(\clE)(m)$ is generated by global sections.
Under this assumption we get a surjection
\begin{eqnarray}
\xymatrix{
H:= H^0(\pi_*\mathcal{E}nd(\clE)(m))\otimes\clE(-m) \ar@{->>}[r] & \clE.
}
\end{eqnarray}

To any pair $(\clF,\phi)$ we can associate a commutative diagram
\begin{eqnarray}
 \xymatrix{
H\times \Gamma \otimes\clE(-m)\ar[r]^-{\tilde{a}}\ar[d]_{\tilde{ev}:=\pi^*ev\otimes\clE(-m)} & V\otimes\clE(-m)\ar[d]^{q}\\
\Gamma \otimes \clE \ar[r]_{\phi} & \clF
}
\end{eqnarray}
where $\tilde{a}:=a\otimes \clE(-m)$ and $a\in \mbox{Hom}(H\times\Gamma,V)$.
In this section we will show that we can use  as  parameter space a suitable locus of $N\times\tilde{Q}$,
where $N$ is the  projectivization of the vector  space  $\mbox{Hom}(H\times\Gamma,V)$.
We start by recalling some useful facts.
\subsection{Quot schemes on projective stacks}
Let $\clX\stackrel{\pi}{\to} X \stackrel{f}{\to} S$ be  a projective  stack.
Let $\tilde{\clQ}$ denote $\mbox{Quot}_{\clX/S}(V\otimes\clE(-m),P)$ and let $Q$ denote
$\mbox{Quot}_{X/S}(F_\clE(V\otimes\clE(-m)),P)$.
By \cite{OlSt03} Proposition 6.2 and \cite{Nir08-Mod} 4.20 there is a closed embedding 
$\iota:\tilde{\clQ}\to Q$.  We know that for any $l\in\bbN$ big enough there is a closed embedding into the Grassmannian
\begin{eqnarray}
j_l: \mbox{Quot}_{X/S}(F_\clE(V\otimes\clE(-m)),P)\hookrightarrow Grass(f_*F_\clE(V\otimes\clE(l-m)),P(l))
\end{eqnarray}
given by the very ample line bundles $det f_{Q*}\clU(l)$, where $\clU$ is the universal quotient sheaf over $Q$
 and $f_{Q}$ is defined as in    the following cartesian diagram
\begin{eqnarray}\label{notation-diag-1}
 \xymatrix{
\clX_{\tilde{\clQ}}\ar[d]_{\pi_{\tilde{\clQ}}}\ar[rrrr] & & & &   \clX \ar[d]^{\pi}\\
X_{\tilde{\clQ}}\ar[d]_{f_{\tilde{\clQ}}}\ar[r]^{\tilde{\iota}}& X_\clQ\ar[d]_{f_Q}\ar[r] & X_{G}\ar[r]\ar[d] & X_{\bbP}\ar[d]\ar[r] & X\ar[d]^{f}\\
\tilde{\clQ}\ar[r]_{\iota} & \clQ\ar[r]_{j_l} & G \ar[r]_{k_l} &  \bbP\ar[r] & S
}
\end{eqnarray}
where $G=Grass(f_*F_\clE(V\otimes\clE(l-m)),P(l))$,  $\bbP=\bbP(\wedge^{P(l)}(f_*F_\clE(V\otimes\clE(l-m)))$,
$j_l$ is the Pl\"ucker embedding and $k_l$ is the Grothendieck embedding.
A result in \cite{Nir08-Mod} identifies  a (relatively) very ample line bundle 
giving an equivariant embedding into the projective space $\bbP(\wedge^{P(l)}(f_*F_\clE(V\otimes\clE(l-m)))$.
\begin{prop}[\cite{Nir08-Mod} Prop 6.2]\label{very-ample-l-b-lem}
The class of invertible sheaves
\begin{eqnarray}
 L_l:=det(\ f_{\tilde{Q}_*}(F_\clE(\tilde{\clU})(l)))
\end{eqnarray}
is very ample for $l$ big enough, where we refer for notation to diagram (\ref{notation-diag-1}) 
and $\tilde{\clU}$ is the universal quotient bundle over $\tilde{\clQ}$.
\end{prop}
As in \cite{HuyLe2} pag. 101, $\tilde{\clU}$ has a natural $GL(V)$-linearization induced by the universal
automorphism of $GL(V)$. As observed in \cite{Nir08-Mod} Lemma 6.3 there is an induced linearization of  $L_l$ 
as its formation commutes  with arbitrary base change.

\subsection{Identification of the parameter space}
Let $N:=\bbP(Hom(H\times \Gamma, V)^\vee) $ be the projective space of morphisms $H\times\Gamma\to V$,
which is polarized by $\clO_N(1)$. 
\begin{lem}
 Let $(\clF,\phi)$ be  a  $\delta$ semistable pair over $\clX$ with $P_\clE(\clF)=P$.
Then it determines a pair $(a,q)$ in $N\times \tilde{Q}$ such that $$q\circ a=\phi\circ \tilde{ev}$$
and  such that $H^0(F_\clE(q(m))\circ \varphi_\clE(V\otimes\clO_X))$ is an isomorphism, where
$\varphi_\clE(V\otimes\clO_X):V\otimes\clO_\clX\hookrightarrow V\otimes\clE nd(\clE)$ is the 
multiplication by the identity endomorphism.
\end{lem}
\begin{pf}
Since $\clF$ is bounded, then for $m$ large enough there is a surjection
\begin{eqnarray}
\pi^* H^0(F_\clE(\clF)(m))\otimes\clE\xrightarrow{G_\clE(ev)} \pi^* F_\clE(\clF)(m)\otimes\clE \xrightarrow{\theta_\clE(\clF)} \clF(m)
\end{eqnarray}
Moreover, $P=h^0(F_\clE(\clF)(m))$. Hence, by choosing an  isomorphism $V\isomto H^0(F_\clE(\clF)(m))$
and by tensoring by $\clO_\clX(-m)$ 
we get a quotient
\begin{eqnarray}
 V\otimes\clE(-m) \to \clF
\end{eqnarray}
namely an element of $\tilde{\clQ}$. Let us consider $\phi:\clE\otimes\Gamma\to \clF$.
By twisting by $\pi^*\clO_X(-m)$,  by applying the functor $F_\clE$ and by taking global sections we get
\begin{eqnarray}
H^0( F_\clE(\clE)(m))\times \Gamma \to H^0(F_\clE(\clF)(m))
\end{eqnarray}
By choosing again an isomorphism $V\isomto H^0(F_\clE(\clF)(m))$ we get  a morphism
$a:H\times\Gamma\to V$. Note that $a$ and $\lambda a$ , $\lambda\in\bbC^*$, come from  isomorphic
pairs. Indeed for any $\lambda\in\bbC^*$
$(\clF,\phi)\simeq (\clF,\lambda\clF)$ by definition.
We observe moreover
$H^0(F_\clE(q(m))\circ\varphi_\clE(V\otimes\clO_X))$ is an isomorphism by construction. The same is true for the relation $q\circ a=\phi\circ \tilde{ev}$.
\end{pf}
We characterize now the locus in $N\times\tilde{Q}$ containing pairs $(a,q)$  yielding
a pair $(\clF,\phi)$.
\begin{prop}
 There is a closed  subscheme $\clW\subseteq N\times\tilde{Q}$ with the following  property:
given a pair $(a,q)\in N\times\tilde{Q}$ 
the map $q\circ \tilde{a}$ factors through $\tilde{ev}$ iff $(a,q)\in\clW$. 
\end{prop}
\begin{pf}
 Same as \cite{Wand10} Prop. 3.4.
\end{pf}
\begin{defn}
We define  $\clZ$ to be the closure in $\clW$  of the open locus of points $(a,q)$
such that $q(V\otimes\clE(-m))$ is pure.
\end{defn}

\section{GIT Construction}
We come now to the GIT construction of the moduli space of $\delta$ semistable pairs.
We observe that $\bbC^*\subset GL(V,k)$ acts trivially on both $N$ and $\tilde{Q}$.
As far as the GIT problem is concerned we can consider the action of the group $PGL(V,k)$ or 
$SL(V,k)$ ($SL(V)$ from now on). Indeed $PGL(V)$ is a quotient of $SL(V)$ by a finite subgroup.
As a consequence, up to taking finite tensor powers, the line bundles linearized for the
actions of the two groups are the same.
There is a natural action of the  group $SL(V)$  on $\clZ$ defined  as follows
\begin{eqnarray}
 g\cdot (a,q)\mapsto ( g^{-1}\cdot a, q\cdot g).
\end{eqnarray}
The line bundles $L_l$ of Lemma \ref{very-ample-l-b-lem} and $\clO_N(1)$  have a natural $SL(V)$ linearization.
We choose as linearized line bundle for the GIT construction
\begin{eqnarray}\label{very-ample-linearized-l.b.}
\clO_\clZ(n_1,n_2):=\pi_{\tilde{Q}}^*L_l^{n_1}\otimes\pi_N^*\clO_N(1)^{n_2}|_{\clZ}.
\end{eqnarray}
We choose $n_1$ and $n_2$ as
\begin{eqnarray}\label{n_1-n_2-choice-eqn}
 \frac{n_1}{n_2}=\frac{P_\clE(l)\delta(m)-P_\clE(m)\delta(l)}{P_\clE(m)+\delta(m)}. 
\end{eqnarray}

\begin{defn}\label{clR-defn}
 Let $\clR\subseteq \clZ$ be the subset of points corresponding to  $\delta$ semistable pairs and such that
 $H^0(F_\clE(q(m))\circ \iota_\clE(V\otimes\clO_\clX))$ is an isomorphism.
\end{defn}
\begin{lem}
 The subset $\clR$ of Definiton \ref{clR-defn} is open and $SL(V)$ invariant. Moreover there is an open subset
$\clR^s\subseteq \clR$  corresponding to $\delta$ stable pairs.
\end{lem}
\begin{pf}
 The proof is almost the same as  \cite{Wand10} {\bf Definition/Lemma 3.5}. 
The proof relies on the fact that the set of Hilbert polynomials 
of purely dimensional quotients destabilizing the pair corresponding to some $(a,q)\in\clZ$ is finite.
To prove this for projective stacks one need two ingredients. The first is 
the Grothendieck Lemma for stacks.  The  second 
is the  fact that given a family of projective stacks and a coherent sheaf over it there exists
 a  finite stratification 
such that the restriction of the given sheaf to each stratum is flat (cfr. \cite{Nir08-Mod} {\bf Proposition 1.13}).
Recall that the Hilbert polynomial over a projective stack is constant in families. 
\end{pf}
\begin{defn}\label{R-bar-defn}
 We define $\overline{\clR}$ as the closure of $\clR$.
\end{defn}

We define the  functor of semistable pairs.
\begin{defn}\label{functor-defn}
Let  the functor
\begin{eqnarray}
 \underline{\clM}_{\clX,\delta}(\clE,P): (Sch/k)^o\to (Sets)
\end{eqnarray}
be defined as follows. For any $S$ in $(Sch/k)$    $\underline{\clM}_{\clX,\delta}(\clE,P)(S)$ is the set of   isomorphism classes of families  $(\clF,\phi)$ of $\delta$ semistable pairs parametrized by $S$   as defined in {\em Definition \ref{ss-pair-family-defn}}, 
such that for  any closed point $s$ of $S$ $(\clF(s), \phi|_{\pi_\clX^*\clE(s)})$ has Hilbert polynomial $P$.
For any $f:S'\to S$ $$\underline{\clM}_{\clX,\delta}(\clE,P)(f):\underline{\clM}_{\clX,\delta}(\clE,P)(S)\to \underline{\clM}_{\clX,\delta}(\clE,P)(S')$$
takes a family $(\clF,\phi)$ over $S$ to $((f\times{\bf 1}_\clX)^*\clF, (f\times{\bf 1}_\clX)^*\phi)$ over $S'$.
 We define
$\underline{\clM}^s_{\clX,\delta}(\clE,P)$ as the subfunctor parametrizing  $\delta$ stable pairs.
\end{defn}

\begin{thm}
 Let $(\clX,\clO_\clX(1),\clE)$ be a polarized smooth  projective   stack. Then the functor $\underline{\clM}_{\clX,\delta}(\clE,P)$ is isomorphic to $[\clR^{ss}/GL(V)]$, where $\clR^{ss}\subseteq \overline{\clR}$ is the subset of GIT semistable points.
\end{thm}
\begin{pf}
The proof uses a quite standard machinery (cfr. \cite{HuyLe2} {\bf Lemma 4.3.1}, \cite{Diac08} {\bf Proposition 3.9}, \cite{Shesh10} {\bf Theorem 4.0.7}).
We sketch it.
We will show that  there is an invertible functor of  categories fibered in groupoids
 between $[\clR^{ss}/GL(V)]$ and $\underline{\clM}_{\clX,\delta}(\clE,P)$. Let us show that $\xi$ exists.
An object of $[\clR^{ss}/GL(V)]$ over $S$ is a diagram
\begin{eqnarray}
 \xymatrix{
P\ar[r]\ar[d] & \clR^{ss}\\
S &
}
\end{eqnarray}
where the horizontal arrow is $GL(V)$ equivariant.
By pulling back the universal family over $\overline{R}\times\clX$ we get a $\delta$ semistable pair
$\phi:\clE\otimes\Gamma\to \clF$ over $P\times\clX$ with an isomorphism
$V\isomto H^0(F_\clE(\clF(m))(s))$ for any closed point $s$ of $S$. The pair comes 
with a $GL(V)$ linearization.  
Therefore both sheaves and the morphism between them descend to $S$. 
The same is true for a morphism between $\delta$ semistable pairs over $P\times\clX$.
We show that there exists a functor $\eta$ in the opposite direction.
Let us draw a diagram to fix the  notations:
\begin{eqnarray}\label{family-over-P-diag}
 \xymatrix{
P\times \clX \ar[r]^{p\times{\bf 1}_\clX}\ar[d]_{\pi_P} & S\times \clX\ar[d]^{\pi_S}\\
P\ar[r]_{p} & S
}.
\end{eqnarray}
Let us consider a $\delta$ semistable pair over $\clX\times S$
$\phi:\clE\otimes\Gamma\to \clF$. Since $\clE$ and  $\clF$ are flat over $S$ by definition of family,
$\pi_{S*}\clE nd(m)$ and $\pi_{S*}\clF\otimes\clE^\vee(m)$ are locally free $\clO_S$-modules and the  morphisms
\begin{eqnarray}
 \pi_{S}^*\mathcal{A}\otimes\clE:=\pi_{S}^*\pi_{S*}(\Gamma\otimes\clE nd(\clE)(m))\otimes\clE\to \clE(m)
\end{eqnarray}
and
\begin{eqnarray}
\pi_{S}^* V\mathcal{B}\otimes\clE:=\pi_{S}^*\pi_{S*}(\clF\otimes\clE^\vee(m))\otimes\clE\to \clF(m)
\end{eqnarray}
are surjective.
Note that $\clA\simeq H\times\Gamma\otimes \clO_S$, $H:=H^0(\clX,\clE nd(\clE)(m))$
Let  $P:=\mathbb{I}som(V,\clB)$ be the frame bundle of $\clB$ (see \cite{HuyLe2} {\bf Example 4.2.3} for the definition). Then  $p:P\to S$  is a $GL(V)$ principal bundle.

 There is a commutative diagram over $P\times\clX$
\begin{eqnarray}
 \xymatrix{
H\times\Gamma\otimes\clE\otimes \clO_{P\times \clX}\ar@{->>}[d] \ar[r] & V\otimes \clE\otimes\clO_{P\times \clX}\ar@{->>}[d]\\
(p\times{\bf 1}_\clX)^* \clE(m)\ar[r] & (p\times{\bf 1}_\clX)^*\clF(m)
}
\end{eqnarray}
where 
\begin{enumerate}
\item[i)] the right vertical arrow corresponds to a morphism    $P\to \tilde{Q}$ and it is obtained   by composing with $(p\times{\bf 1}_\clX)^*\clF(m)$
with $\pi_P^*\mathfrak{h}^{-1}$, where $\mathfrak{h}$ is the universal isomorphism
  over $P$;\\
\item[ii)] the upper  horizontal arrow corresponds to a morphism  $P\to N$ and it is obtained
by post-composing the natural morphism $H\times\Gamma\otimes \clO_{P\times\clX}\to (p\times{\bf 1}_\clX)^*\pi_S^*\clB$
with $\pi_P^*\mathfrak{h}^{-1}$.
\end{enumerate}
By Cohomology and base change  for Deligne-Mumford stacks (see \cite{Nir08-Mod} {\bf Theorem 1.7}) there are isomorphisms
$H\times\Gamma\otimes\clO_{P} \isomto p^*\clA$
and $V\otimes \clO_{P} \isomto p^*\clB$.
Let $p$ be  a closed point of $P$. 
 Then
\begin{eqnarray}\label{isom-cohom-prop}
H\times\Gamma\simeq  p^*\clA|_p\simeq H^0(F_\clE(\clE\otimes\Gamma(m))(p))\nonumber\\
V\simeq  p^*\clB|_p\simeq H^0(F_\clE(\clF(m))(p))
\end{eqnarray}
Then  the family of diagram  (\ref{family-over-P-diag}) with the  property (\ref{isom-cohom-prop})
and its natural $GL(V)$ linearization
 provides 
a $GL(V)$-equivariant morphism to $\clR^{ss}$.
Such a morphism descends to a morphism $S\to [\clR^{ss}/GL(V)]$.
In a similar way we can reconstruct morphisms.
Let  $f:S\to T$ be  a morphism, 
let $(\clF,\phi)$   and $(\clF',\phi')$ be families of pairs over $S\times\clX$ and $T\times \clX$,
$\tau:(\clF,\phi)\to (p\times{\bf 1}_\clX)^*(\clF',\phi')$ be a morphism over $S\times\clX$. 
By cohomology and base change along the diagram
\begin{eqnarray}
 \xymatrix{
S\times \clX \ar[r]\ar[d] & T\times\clX\ar[d]\\
S\ar[r] & T
}
\end{eqnarray}
 we see that the families produce $GL(V)$ principal bundles $P$ and $P'$ such that
$P\isomto f^*P'$. 
By cohomology and base change along the diagram
\begin{eqnarray}
 \xymatrix{
f^*P'\times\clX\ar[r]^u\ar[d] & P'\times\clX\ar[d]\\
f^*P'\ar[r] & P'
}
\end{eqnarray}
 we see that there is a canonical  isomorphism between the pair over $P$  and the pullback
of the pair over $P'$. Then we get two isomorphic $GL(V)$ equivariant maps 
$g:P\to \clR^{ss}$ and $g':f^*P'\to  \clR^{ss}$ such that
$g'\circ u=g$.
We conclude that the functor $\eta$ is defined.
It is straightforward  to check that $\xi$ and $\eta$ are the inverse of each other.  
 \end{pf}

\begin{lem}\label{categorical-quotient-corepr-lem}
 Let $M$ be a scheme which is  a categorical quotient for the $SL(V)$ action on $\clR$. Then it corepresent
the functor $\underline{\clM}_{\clX,\delta}(\clE,P)$.
\end{lem}
\begin{pf}
 The functors  $\underline{\clM}_{\clX,\delta}(\clE,P)$  and $[\clR/GL(V)]$
are isomorphic. Hence to corepresent either functor is the same. 
\end{pf}

We state the theorem relating GIT semistability of points of the parameter space to
$\delta$ semistability of the corresponding pairs. The proof will be given by {\em Theorems \ref{GIT-implies-dss-thm} and 
\ref{delta-ss-implies-GIT-ss-thm}}. 
\begin{thm}\label{main-thm}
For $l$ large enough
 the subset of points in the  closure $\overline{\clR}$ of  $\clR$ which are semistable with respect
to $\clO_{\clZ}(n_1,n_2)$ and its $SL(V)$ linearization coincides with the 
subset of points corresponding to $\delta$ semistable pairs.
\end{thm}
\begin{thm}
 Let $(\clX,\clO_\clX(1),\clE)$ be a polarized smooth  Deligne-Mumford  stack. 
Then there exists a moduli space    $M_{\clX,\delta}(\clE,P)$
 of $\delta$ semistable pairs. Two pairs correspond to the same point in  $M_{\clX,\delta}(\clE,P)$ if and only if they are S-equivalent. Moreover there is an open
subset $M_{\clX,\delta}^s(\clE,P)$ corresponding to stable pairs. It is a fine moduli space for $\delta$ semistable pairs.
\end{thm}
\begin{pf}
The prrof is an in \cite{Wand10} {\bf Theorem 3.8}.
 We proved in  {\em Lemma \ref{categorical-quotient-corepr-lem}}
 that a categorical quotient of $\cl$ for the $SL(V)$ action
corepresent the functor $\underline{\clM}_{\clX,\delta}(\clE,P)$.
General GIT theory and  {\em Theorem \ref{main-thm}} 
provide a categorical quotient of $\clR$ which in particular is a good quotient.
By using arguments analogous to \cite{HuyLe95} {\bf Proposition 3.3} or \cite{HuyLe2} {\bf Theorem 4.3.3}
one proves that closure of the orbit of a $\delta$ semistable pair $(\clF,\phi)$  contains 
the pair $(gr^{JH}(\clF), gr^{JH}(\phi))$. Moreover the orbit of  $(\clF,\phi)$ is closed if and only if
 $(\clF,\phi)$ is polystable. In proving the last fact a semicontinuity result in used, which also holds
for Deligne-Mumford stacks (cfr. \cite{Nir08-Mod}) {\bf Theorem 1.8}. The second statement follows then
from properties of good quotients.
Moreover there is a universal family over $\clX\times M^s_{\clX,\delta}(\clE,P)$, which is implied 
by the existence of  a line bundle over $\clR^s$ of weight $1$ for the action of $\bbG_{m}\subseteq GL(V)$
(cfr. \cite{HuyLe2} {\bf Section 4.6}). Such a line bundle is provided by $\clO_N(1)$.
\end{pf}

\subsection{GIT computations}
We give in the following proposition numerical conditions implied by GIT semistability.
We recall some useful standard results.
The action of a  1-parameter subgroup of $SL(V)$ $\lambda$ can be diagonalized. Hence  $V$ splits in eigenspaces
for the eigenvectors of $\lambda$. The subspaces  $V_{\leq n} =\oplus_{i=1}^n V_i\subseteq V$ give an
ascending filtration.
Let $q:V\otimes \clE(-m)\to \clF$ a closed point of $\tilde{Q}$.
We get a corresponding filtration  $\clF_{\leq n}$ of $\clF$, where $\clF_{\leq n}=q(V_{\leq n}\otimes\clE(-m))$.
Let $\clF_n=\clF_{\leq n}/\clF_{\leq n-1}$ the graduate pieces.
Let
\begin{eqnarray}
 \overline{q}: V\otimes\clE(-m) \to \oplus_{i=1}^n \clF_{\leq n}.
\end{eqnarray}
be a closed point of $\tilde{Q}$.
Then we have the following result generalizing \cite{HuyLe2} Lemma 4.4.3 and proven in \cite{Nir08-Mod} Lemma 6.11.
\begin{lem}
 The quotient $[\overline{q}]$ is  $\lim_{t\to 0}\lambda(t)q$ in the sense of the Hilbert-Mumford criterion.
\end{lem}
As in the case of scheme the action of $\lambda$ on the fiber of $L_l$ is characterized as follows.
\begin{lem}[{\bf \cite{Nir08-Mod} Lem. 6.12}]\label{quot-weight}
 The action of $\bbG_{m,k}$ via the representation $\lambda$ on the fiber of $L_l$ at the point $[\overline{q}]$
is given by the weight
\begin{eqnarray}
 \sum_n P_\clE(\clF_n(l)).
\end{eqnarray}
\end{lem}

We fix some notation. Let $W:=H^0(\clO_X(l-m))$. We denote by $q'$ the following morphism induced by $q$:
\begin{eqnarray}
 q':V\otimes W\xrightarrow{H^0(\varphi_\clE(l))\circ H^0(\phi_\clE(l))} H^0(\clF\otimes\clE^\vee(l))
\end{eqnarray}
and by $q''$ its $r$-th antisymmetric product
\begin{eqnarray}\label{q''-eqn}
 q'':\wedge^r(V\otimes W) \to det\ H^0(\clF\otimes\clE^\vee(l)).
\end{eqnarray}

\begin{prop}
 Let $(a,q)$ be a point in $\overline{\clR}$. For $l$ large enough $(a,q)$ is GIT semistable with respect
to $\clO_{\clZ}(n_1,n_2)$ if and only if the following holds. For any non trivial subspace 
$U\subseteq V$ we have
\begin{eqnarray}\label{GIT-num-cond-1-eqn}
 dim\ U[n_1 P(l) - n_2] (\leq) P(m)[n_1 dim\ q'(U\otimes W)-\epsilon(U)n_2]
\end{eqnarray}
where $\epsilon(U)=1$ if $im\ a\subseteq U$ and 0 otherwise.
\end{prop}
\begin{pf}
We use Hilbert-Mumford criterion.
Let $\lambda:\bbC^*\to SL(V)$ be a one parameter subgroup. Let us choose  a basis of $V$ $v_1,..,v_p$
such that $v_i\cdot\lambda(t)=t^{\gamma_i}v_i$, $\gamma_{i+1}\geq \gamma_i$, $i=1..,p-1$.
Recall that $\lambda$ is completely specified by a weight vector $(\gamma_1,...,\gamma_p)\in\bbZ^p$
 Moreover $\sum_{i=1}^p \gamma_i=0$. 
The number $\mu(q,\lambda)$ is computed by Lemma \ref{quot-weight}.
Let us compute $\mu(a,\lambda)$, $a\in N=Proj(Sym^{\bullet} Hom(H\times\Gamma,V)^\vee)$.
Sections of $\clO_N(1)$ are $Hom(H\times\Gamma,V)^\vee$.
Let $\xi:H\times\Gamma\to V$ be a homomorphism. We can write it in matricial form
as $\xi=\sum_{i,j}\alpha_{ji} w_j^\vee\otimes v_i$, $w_j\in H\times\Gamma$.
Then
\begin{eqnarray}
 \mu(a,\lambda)=\mbox{max}\{\gamma_i|\alpha_{ij}\neq 0\}.
\end{eqnarray}
Equivalently  $\mu(a,\lambda)=\gamma_i$ where $i=\mbox{min}\{\ i|\mbox{im}\ a\subseteq \langle v_1,.,v_i\rangle  \}$.
By the Hilbert-Mumford criterion (semi) stability of $(a,q)$ requires that
\begin{eqnarray}
 \mu(q'',\lambda)n_1+\mu(a,\lambda)n_2(\geq)0\quad.
\end{eqnarray}
Given a base $v_1,..,v_p$ of $V$, let consider weight vectors of the form
\begin{eqnarray}\label{sample-weight-vector}
 \gamma^{(i)}=(\underbrace{i-p,..,i-p}_{i},\underbrace{i,..,i}_{p-i})
\end{eqnarray}
 Any other weight vector can be expressed as a finite non negative linear combination of weight vectors in this class. In fact 
\begin{displaymath}
 \gamma_k=\sum_{i=1}^{p-1} \left( \frac{\gamma_{i+1}-\gamma_i}{p}\right) \gamma^{(i)}_k.
\end{displaymath}
Let us start by evaluating $\mu(q,\lambda)$. By Lemma \ref{quot-weight} we get
\begin{eqnarray} 
 \mu(q,\lambda)=\sum_n n P(\clF_n(l))=(i-p)\sum_{n\leq i} P(\clF_n(l)) + i \sum_{i<n\leq p} P(\clF_n(l))=\nonumber\\
i\rho - p\psi(i),
\end{eqnarray}
where $\rho=h^0(\clF\otimes\clE^\vee(l))$ and  $\psi(i)= \sum_{n\leq i} P(\clF_n(l))=\mbox{dim}(q'(\langle v_1,..,v_i\rangle\otimes W))$.
Here the second equality holds because the graded pieces are also bounded. 
Let us come to $\mu(a,\lambda)$.
Note that  $\mu(a,\lambda)=i-p$ if $\mbox{im}\ a\subseteq \langle v_1,..,v_i\rangle$,
and that $\mu(a,\lambda)=i$ otherwise. More compactly,  $\mu(a,\lambda)=i-\epsilon(i)p$,
where $\epsilon(i)=1$ if $\mbox{im}\ a\subseteq \langle v_1,..,v_i\rangle$, 0 otherwise.
Summing up the Hilbert-Mumford criterion implies
\begin{eqnarray} \label{1st-inequality}
 i(\rho n_1 -n_2) (\leq ) p(\psi(i)n_1 - \epsilon(i)n_2)\quad 
\end{eqnarray}
We observe that the above equation does neither depend on the base of $V$ we chosed nor
on the particular weight vector we evaluated the Hilbert-Mumford criterion on.
Let $U\subseteq V$ the vector subspace generated by $v_1,.,v_i$. Since $\rho=P_\clE(l)$
and $\mbox{dim} V= P_\clE(m)$  the inequality (\ref{1st-inequality}) can be rewritten as
\begin{eqnarray}\label{GIT-num-cond-1-eqn}
 \mbox{dim}\ U(P_\clE(l)n_1 -n_2) (\leq) \mbox{dim}\ V(\mbox{dim}(q'(U\otimes W)n_1 - \epsilon(U)n_2)\quad 
\end{eqnarray}
where $\epsilon(U)=1$ if $\mbox{im}\ a\subseteq U$ and 0 otherwise.
\end{pf}
\begin{notation}\label{FU-notation}
 Let $U\subseteq V$  be a sub vector space. We denote by $\clF_U$ the subsheaf of $\clF$
generated by $q(G_\clE(U))$.
\end{notation}
\begin{cor}\label{injectivity-cor}
 For  any GIT-semistable point  $(a,q)$ the induced morphism 
$$q(m)\circ\iota_\clE(V\otimes\clO_\clX): V\otimes\clO_\clX\to V \otimes\mathcal{E}nd\ \clE\to \clF\otimes\clE^\vee(m)$$
 is injective.
In particular $\mbox{dim}(V\cap H^0(F_\clE(\clF)(m)))\leq h^0(F_\clE(\clF)(m))$. Moreover $q'$ is injective
and for any subspace $U\subseteq V$ $\mbox{dim}\ q'(U\otimes W)\leq h^0(\clF_U(l))$
where $\clF_U=q(U(-m))$.
\end{cor}
\begin{pf}
Let $U\subseteq V$ be the kernel of $q(m)\circ\iota_\clE(V\otimes\clO_\clX)$. Then $\epsilon(U)=0$, otherwise
$\phi$ would be zero. Moreover 
$\mbox{dim}\ q'(U\otimes V)=0$. Substituting in equation 
 (\ref{GIT-num-cond-1-eqn}) we get $\mbox{dim}\ U\leq 0$  because $n_1\rho -n_2\geq 0$
by the choice (\ref{n_1-n_2-choice-eqn}).
\end{pf}
We give another numerical characterization  for GIT semi-stability.
\begin{prop}
 For sufficiently large $l$  a point $(a,q)$ is GIT-semistable  if for any
$U\subseteq V$ the following polynomial equation holds
\begin{eqnarray}\label{GIT-num-cond-2-eqn}
 \mbox{dim}\ U(P_\clE(l)n_1 -n_2) (\leq) \mbox{dim}\ V(\mbox{dim}(P_\clE(\clF_U(l))n_1 - \epsilon(U)n_2)
\end{eqnarray}
\end{prop}
\begin{pf}
 It is enough to prove that equation  (\ref{GIT-num-cond-2-eqn}) implies (\ref{GIT-num-cond-1-eqn}). 
Note that subsheaves of the form $\clF_U$ are bounded. Hence for $l$ large enough
$P_\clE(\clF_U(l))=h^0( \clF_U\otimes\clE^\vee(l))=\mbox{dim}\ q'(U\otimes W)$.  We are left to show  that $\epsilon(U)=1\Leftrightarrow \epsilon(\clF_U)=1$.
For any $U\subseteq V$  $\epsilon(U)=1\Rightarrow \epsilon(\clF_U)=1$. Let us prove the opposite implication.
It can happen  that $\mbox{Im}\ a\subsetneq U$
but  $\mbox{Im}\ \phi\subseteq \clF_U$. Let  $U'$ be the subspace of $V$ generated by $U$
and by $\mbox{Im}\ a$. Then the inequality holds with $\mbox{dim}\ U'$ replacing $\mbox{dim}\ U$,
hence it holds a fortiori for $\mbox{dim}\ U$.
\end{pf}
\begin{prop}
 For sufficiently large $l$  a point $(a,q)$ is GIT-semistable if and only if for any
$U\subseteq V$ the following polynomial equation holds
\begin{eqnarray}\label{GIT-num-cond-3-eqn}
 P(\mbox{dim}\ U + \epsilon(\clF_U)\delta(m))+\delta (\mbox{dim}\ U -\epsilon(\clF_U)P(m))(\leq) P_{\clF_U}(P(m)+\delta(m))
\end{eqnarray}
\end{prop}
\begin{pf}
The proof is as in \cite{Wand10}.
Take $l$ large enough such that inequality (\ref{GIT-num-cond-2-eqn}) holds as an inequality of polynomials.
Then put
 \begin{eqnarray}
  \frac{n_1}{n_2}=\frac{P_\clE(l)\delta(m)-P_\clE(m)\delta(l)}{P(m)+\delta(m)}
\end{eqnarray}
\end{pf}
\begin{thm}\label{GIT-implies-dss-thm}
 For sufficiently large $l$ if  a point $(a,q)$ in $\overline{\clR}$ is GIT-semistable then the corresponding pair
$(\clF,\phi)$ is $\delta$ semi-stable and $H^0(q(m)\circ \iota_\clE(V\otimes\clO_\clX))$ is an isomorphism.
In particular any GIT-semistable point corresponds to a pair with torsion-free sheaf.
\end{thm}
 \begin{pf}
Note that by {\em Corollary \ref{injectivity-cor}} $V\otimes\clO_\clX \to \clF\otimes\clE^\vee(m)$ is injective.
Hence for dimensional reasons $V\to H^0(\clF\otimes\clE^\vee(m))$ is an isomorphism. 
  Let $\clF'\subseteq \clF$ be a subsheaf. Let $U=V\cap H^0(F_\clE(\clF'))$. Then
$\mbox{dim}\ U=h^0(F_\clE(\clF'))$. 
Let $\clF_U\subseteq \clF'$ as in {\em Notation \ref{FU-notation}}. We observe that 
$\epsilon(\clF_U)=1$ iff $\epsilon(\clF')=1$. The ``if'' direction holds because 
$\epsilon(\clF')=1$ implies $\epsilon(U)=1$ since  $U=V\cap H^0(F_\clE(\clF'))$.
This in turn implies $\epsilon(\clF_U)=1$.
By taking the leading coefficients in the polinomial equation (\ref{GIT-num-cond-3-eqn})
we get
\begin{eqnarray}\label{leading-coeff-GIT-ss-eqn}
 \mbox{dim} U +\epsilon(\clF')\delta(m)\leq \frac{r'}{r}(P_\clE(m)+\delta(m))
\end{eqnarray} 
By  {\em Lemma \ref{deform-sheaf-lem}} there exists a morphism 
$\psi:\clF\to\clH$, where $\clH$ is pure, the kernel is a torsion subsheaf $\clT$  and 
$P_\clE(\clF)=P_\clE(\clH)$. 
There is an induced homomorphism $\phi_H:\clE\otimes\Gamma\to \clH$ which is non vanishing.
If it was, then (\ref{leading-coeff-GIT-ss-eqn}) would be violated for $\clF'=\clT$.
 Let $\clH''$ be a quotient of $\clH$, and $\clH'$ the corresponding kernel.
Let $\clF''$ the image of $\clF$ in $\clH''$ and let $\clF'$ be the corresponding kernel. Then 
\begin{eqnarray}\label{clH-is-ss}
 h^0(\clH''\otimes\clE^\vee(m))+\epsilon(\clH'')\delta(m)\geq\nonumber\\
h^0(\clF''\otimes\clE^\vee)+\epsilon(\clF'')\delta(m)\geq\nonumber\\
\mbox{dim} V + \delta(m) - (\mbox{dim} U + \epsilon(\clF'))\nonumber\geq\\
P_\clE(m)+\delta(m) - \frac{r_{\clE,\clF'}}{r_{\clE,\clF}} (P_\clE(m)+\delta(m))\geq\nonumber\\
\frac{r_{\clE,\clF''}}{r_{\clE,\clF}} (P_\clE(m)+\delta(m))=\nonumber\\
\frac{r_{\clE,\clH''}}{r_{\clE,\clF}} (P_\clE(m)+\delta(m))\ \ 
\end{eqnarray}
It follows that $(\clH,\phi_H)$  is $\delta$ semistable by {\em Proposition \ref{num-cond-d-ss-prop}}, therefore belongs to  a bounded family.
Consequently it is $m$-regular.
 Let us consider the sheaf  $\psi(\clF)\subseteq \clH$. It is a quotient of $\clF$. By
(\ref{clH-is-ss}) $h^0(\psi(\clF)\otimes\clE^\vee(m))\geq h^0(\clH\otimes\clE^\vee(m))=P_\clE(m)$,
while the opposite inequality is obvious. Hence $V\otimes\clE^\vee(-m)\to \clH$ is surjective,
which implies $\clF\to \clH$ is surjective.
By semistability of $\clH$  and $P_\clE(\clF)=P_\clE(\clH)$ it follows that $\psi$ is  an isomorphism.
Eventually, it is an obvious consequence that $V\to H^0(\clH\otimes\clE^\vee(m))$ is an isomorphism, because
$V\to H^0(\clF\otimes\clE^\vee(m))$ is.
\end{pf}
We prove now the inverse implication.
\begin{thm}\label{delta-ss-implies-GIT-ss-thm}
 Let $(\clF,\phi)$ be a $\delta$-semistable pair such that $H^0(F_\clE(q(m))\circ\iota_\clE(V\otimes\clO_\clX)):V\to H^0(F_\clE(\clF))$ is an isomorphism. Then the corresponding point  is GIT-semistable.
\end{thm}
\begin{pf}
If $\delta$ stability holds we can prove  (\ref{GIT-num-cond-3-eqn}) by only considering the leading coefficients.
In this case, given $\clF'\subseteq\clF$ we define $U=V\cap H^0(\clF'\otimes\clE^\vee(m))$.
Then $\mbox{dim} U\leq h^0(\clF'\otimes\clE^\vee(m))$.
Recall $\epsilon(\clF')=1\Leftrightarrow \epsilon(\clF_U)$ where $\clF_U$ is the sheaf generated by $U$ via $q(m)\otimes\clE^\vee\circ\varphi_\clE$.
 This inequality and  the $\delta$ stability condition yield
\begin{displaymath}
 \mbox{dim} U + \epsilon(\clF')\delta(m) < \frac{r_{\clE,\clF'}}{r_{\clE,\clF}}(P_\clE(m)+\delta(m)),
\end{displaymath}
which implies (\ref{GIT-num-cond-3-eqn}) as a polynomial inequality.
The rest of the proof taking into account $\delta$ semistability proceeds exactly
as in \cite{Wand10} {\bf Theorem 4.7}.
In particular a pair $(\clF,\phi)$ is GIT stable only if it is $\delta$ stable.
\end{pf}

%\bibliography{../BIBTEX/el_biblio2.bib}
\end{document}